\documentclass[11pt]{article}
\usepackage{epsf}
\usepackage{color}
\usepackage{epsfig}

\usepackage{graphicx}
\usepackage{amsmath}
\usepackage{latexsym}
\usepackage{amssymb}

\usepackage{a4}

\newtheorem{theorem}{Theorem}

\newcommand{\R}{\mathbb{R}}
\newcommand{\eps}{\varepsilon}

\newcommand{\e}{{\rm e}}

\begin{document}

\title{New families of symplectic splitting methods for numerical integration in dynamical astronomy}

\author{
 S. Blanes$^{1}$\thanks{Email: \texttt{serblaza@imm.upv.es}}
   \and
 F. Casas$^{2}$\thanks{Email: \texttt{Fernando.Casas@uji.es}}
   \and
 A. Farr\'es$^{3}$\thanks{Email: \texttt{afarres@imcce.fr}}
   \and
 J. Laskar$^{3}$\thanks{Email: \texttt{laskar@imcce.fr}}
   \and
 J. Makazaga$^{4}$\thanks{Email: \texttt{Joseba.Makazaga@ehu.es}}
   \and
 A. Murua$^{4}$\thanks{Email: \texttt{Ander.Murua@ehu.es}}
   }

\date{}
\maketitle

\begin{abstract}

We present new splitting methods  designed for the numerical integration
of near-integrable Hamiltonian systems, and in particular for planetary N-body problems,
when one is interested in very
accurate results over a large time span. We derive in a systematic way an independent set of
necessary and sufficient conditions to be satisfied by the
coefficients of splitting methods to achieve a prescribed order of
accuracy. Splitting methods satisfying such (generalized) order
conditions are appropriate in particular for the numerical simulation
of the Solar System described in Jacobi coordinates. We show that, when using Poincar\'e Heliocentric
coordinates,  the same order of accuracy may be obtained by imposing
an additional polynomial equation on the coefficients of the splitting
method. We construct several splitting methods appropriate for each of the two
sets of coordinates by solving the corresponding
systems of polynomial equations and finding the optimal solutions. The
experiments reported here indicate that the efficiency of our new schemes
is clearly
superior to previous integrators when high accuracy is required.

\vspace*{0.6cm}

\begin{description}
 \item $^1$Instituto de Matem\'atica Multidisciplinar,
  Universitat Polit\`ecnica de Val\`encia, E-46022  Valencia, Spain.
 \item $^2$Institut de Matem\`atiques i Aplicacions de Castell\'o and 
   Departament de Matem\`atiques, Universitat Jaume I,
  E-12071 Castell\'on, Spain.
 \item $^3$Astronomie et Syst\`emes Dynamiques, IMCCE-CNRS UMR8028, Observatoire de Paris, UPMC,
 77 Av. Denfert-Rochereau, 75014 Paris, France.
 \item $^4$Konputazio Zientziak eta A.A. saila, Informatika
Fakultatea, EHU/UPV, E-20018, Donostia/San Sebasti\'an, Spain.
\end{description}

\end{abstract}

\section{Introduction}    \label{Introduction}

Symplectic integrators have several features that turn out to be 
particularly appropriate when integrating numerically for long times 
evolution problems in dynamical astronomy. 
They preserve by construction the symplectic structure of the original Hamiltonian problem, so that the 
numerical solution inherits the qualitative properties of the exact one \cite{sanz-serna94nhp}.  In particular, by using backward error
analysis, it is possible to prove that this
numerical solution is in fact exponentially close to the exact solution of a modified Hamiltonian.  Moreover,
although the
energy is not conserved along the trajectory, the error introduced by a symplectic method of order $r$ used with
constant step size $\tau$ is of order $\mathcal{O}(\tau^r)$ for exponentially long time intervals under
rather general assumptions, whereas the error in
position typically grows linearly with time \cite{hairer06gni}.

Assume that, as is often the case, the Hamiltonian function is of the form $H(q,p) = T(p) + U(q)$, where
the potential energy $U(q)$ depends on positions and the kinetic energy $T(p)$ is a function of the
conjugate momenta. Then the equations of motion corresponding to $T(p)$ are trivially solvable, and the same happens
with $U(q)$. By composing the flows of these two special Hamiltonian systems one gets a symplectic first
order approximation to the exact flow. This simple composition constitutes an example of a symplectic splitting method.
Higher order approximations can be obtained by composing the flows corresponding to $T(p)$ and $U(q)$ 
with certain coefficients 
obtained by solving the so-called order conditions \cite{mclachlan02sm}. 
There exist in the literature a vast number of high order 
integrators constructed along this line (see, e.g., \cite{blanes08sac}, \cite{hairer06gni},  and references therein).

The non-relativistic gravitational N-body problem, in particular, belongs to this class of systems. If one considers the
motion of $n+1$ particles (the Sun, with mass $m_0$, and $n$ planets with masses $m_i$, $i=1,\ldots,n$) only
affected by their mutual gravitational interaction, the corresponding equations of motion can be derived from the
Hamiltonian
\begin{equation}   \label{n-body}
  H = \frac{1}{2} \sum_{i=0}^n \frac{\|\mathbf{p}_i\|^2}{m_i} - G \sum_{0 \le i < j \le n} \frac{m_i m_j}{\|\mathbf{q}_i - 
  \mathbf{q}_j\|},
\end{equation}
where $\mathbf{q}_i$ and $\mathbf{p}_i = m_i \, \mathbf{\dot{q}}_i$ 
denote the position and momenta of the $n+1$ bodies in
a barycentric reference frame.  Typically, the planets evolve around the central mass following almost Keplerian
orbits, so that by an appropriate change of coordinates one can rewrite the Hamiltonian (\ref{n-body}) as 
$H = H_K + H_I$, where in some sense $|H_I| \ll |H_K|$, or equivalently, 
as the sum of the Keplerian motion of each planet around the central mass and a small perturbation due to the
gravitational interaction between planets. Jacobi and Heliocentric coordinates constitute paradigmatic examples
of canonical set of coordinates possessing this feature. Thus, the Hamiltonian (\ref{n-body}) written as $H = H_K + H_I$
is a particular example of a near-integrable Hamiltonian system, i.e, it can be expressed as 
\begin{equation}    \label{intro.1}
   H(q,p;\eps) = H^{[a]}(q,p) + \eps H^{[b]}(q,p),
\end{equation}   
where $\eps \ll 1$ and $H^{[a]}$ is exactly integrable. It makes sense, then, to take into account this special structure
when designing integration methods to approximate its dynamics. The idea consists in constructing 
splitting schemes as compositions of the flows corresponding to $H^{[a]}(q,p)$ and $H^{[b]}(q,p)$, assuming that
they are explicitly computable or sufficiently well approximated \cite{kinoshita91sia,wisdom91smf}.  
In fact, 
since the parameter $\eps$ is small, it is possible to design
methods which behave in practice as high order integrators with less severe restrictions concerning the
order conditions than the usual split into kinetic and potential energy. 
This approach was systematically pursued by McLachlan \cite{mclachlan95cmi}, obtaining families of
splitting schemes of order 2 and 4 which eliminate the most relevant error terms in $\eps$, and 
further analyzed by
Laskar \& Robutel \cite{laskar01hos} in the context of planetary motion.

By incorporating the idea of processing, even more efficient schemes can be constructed for the Hamiltonian
(\ref{intro.1}) \cite{blanes00psm}. In that case, both
the kernel and the processor are taken as compositions of the flows associated 
with $H^{[a]}$ and $H^{[b]}$, so that the exactly symplectic character of the integration scheme is ensured.
With this approach, all terms of first order in $\eps$ in the truncation
error expansion can be annihilated with the processor \cite{mclachlan96mos,wisdom96sco}.

Although the symplectic methods developed in \cite{laskar01hos} and \cite{mclachlan95cmi}  for near-in\-te\-gra\-ble Hamiltonian 
systems have proved their usefulness in long term integrations of the Solar System \cite{laskar09eoc}, 
the design of new
and more efficient higher order integrators is of interest for numerical 
simulations of its evolution over large  time spans, either by speeding up the algorithms or by providing
better  accuracy in the position of the different objects. 
Relevant examples where the new integrators could be useful include  the numerical
integration of the Solar System for more than 60 million years backward in
time to cover the Palaeogene period to determine insolation quantities of the Earth and calibrate
paleoclimatic data, studies of the planetary orbits over several billion years, etc. \cite{laskar04alt}.
To this purpose, it is essential that the numerical solutions obtained
are not contaminated by error accumulations along the integration and that the computations are
done in a reasonable time.

These long-time numerical integrations can be combined with standard techniques of classical perturbation
theory, such as the expansion of  the equations of motion up to a certain order in the perturbation parameters
and the use of averaging (see, e.g. \cite{laskar86sto}).

The purpose of this work is to present new families of symplectic splitting methods specifically
designed for Hamiltonian systems of the form (\ref{intro.1}) appearing in many problems of
dynamical astronomy, when one is
interested in highly accurate results over a large time span. The schemes
we propose will be useful in particular in  
the long time integration of the Solar System, both in Jacobi and Poincar\'e Heliocentric coordinates, and are 
more efficient than the schemes designed in \cite{laskar01hos} and \cite{mclachlan95cmi}. Although they
involve the computation of more elementary flows
per step than other methods,  their small error terms allow to use larger
steps, which results in more efficient schemes. Obtaining these new methods requires deriving
previously the necessary and sufficient order conditions to be
satisfied by the coefficients (which is done here in a
systematic way) and then solving these
polynomial equations to get the best solutions according with some appropriately chosen optimization
criteria. This is discussed in more detail in sections \ref{sec.2}, \ref{sec.3} and the appendix, whereas in 
section \ref{sec.4} we consider the application
of the new schemes to the integration of the Solar System. The new
methods obtained in section~\ref{sec.3} are suitable to be applied
when using Jacobi coordinates and also Poincar\'e Heliocentric
coordinates.

It is worth stressing that,  while the main motivation of this work is the long time 
integration of Hamiltonian problems arising in dynamical astronomy, and in particular in planetary systems,
the new symplectic splitting methods 
obtained here can also be applied to more general perturbed differential
equations arising in different fields 
when high accuracies are required.

\section{Order conditions}
\label{sec.2}

\subsection{Preliminaries}

To establish the framework for the construction and analysis of the new families of integrators,
we consider a generic differential equation of the form
\begin{equation}     \label{eq:1}
      x^\prime = f^{[a]}(x) + \eps f^{[b]}(x), \qquad \quad
        x(0)=x_0\in\mathbb{R}^D,
\end{equation}
where $|\eps|\ll 1$ and each part
\begin{equation}     \label{eq:2}
      x^\prime = f^{[a]}(x),  \qquad  x^\prime = \eps f^{[b]}(x)
\end{equation}
is exactly solvable (or can be numerically solved up to round off
accuracy) with solutions
\[
  x(\tau)=\varphi^{[a]}_{\tau}(x_0), \qquad\quad  x(\tau)=\varphi^{[b]}_{\tau}(x_0)
\]
respectively, at $t=\tau$, the time step. If we denote by $\varphi_{\tau}(x_0)$ the exact solution of (\ref{eq:1}), it
is well known that $\psi_{\tau} = \varphi^{[b]}_{\tau} \circ \varphi^{[a]}_{\tau}$ provides a first-order approximation, i.e.,
$\psi_{\tau}(x_0) = \varphi_{\tau}(x_0) + \mathcal{O}(\tau^2)$ and that higher order approximations can be obtained by
taking more compositions in $\psi_{\tau}$,
\begin{equation}  \label{eq:3}
    \psi_{\tau} = \varphi^{[a]}_{a_{s+1} \tau} \circ  \varphi^{[b]}_{b_{s} \tau}\circ \varphi^{[a]}_{a_{s} \tau}\circ
 \cdots\circ  \varphi^{[b]}_{b_{1}\tau} \circ \varphi^{[a]}_{a_{1}\tau}
\end{equation}
for appropriately chosen coefficients $a_i,b_i$.
The splitting method $\psi_{\tau}$ is said to be of order $r$ if for all $x \in \mathbb{R}^D$,
\begin{equation}
  \label{eq:phipsi}
 \psi_{\tau}(x) = \varphi_{\tau}(x) +
\mathcal{O}(\tau^{r+1}) \quad \mbox{as} \quad  \tau \rightarrow 0.
\end{equation}
It is straightforward to check that the method is at least of order 1 for arbitrary problems of the form (\ref{eq:1}) if and only if the coefficients $a_i,b_i$ satisfy the {\em consistency condition}
\begin{equation}
\label{eq:consistency}
  \sum_{i=1}^{s+1}a_i=1, \qquad\quad   \sum_{i=1}^{s}b_i=1.
\end{equation}

We are mainly interested in \emph{symmetric} methods, that is, integrators verifying
$\psi_{-\tau}=\psi_{\tau}^{-1}$, or equivalently $a_{s+2-i} = a_i$, $b_{s+1-i} = b_i$ (so that the
composition (\ref{eq:3}) is left-right palindromic). In that case, they are automatically of even order.  In particular, if a symmetric method satisfies the consistency condition (\ref{eq:consistency}), then it is at least of order 2.

 Since the last flow
$\varphi^{[a]}_{a_{s+1} \tau}$ can be concatenated with the first $\varphi^{[a]}_{a_{1} \tau}$ at the next step when
scheme (\ref{eq:3}) is iterated, the number of flows $\varphi^{[a]}_{\tau}$ and $\varphi^{[b]}_{\tau}$ per step is
precisely $s$. This number is usually referred to as the number of stages in the composition.

\subsection{Deriving the order conditions via the BCH formula}
\label{ss:BCH}
The conditions  that the coefficients $a_i, b_i$ must satisfy for a splitting method to be of order $r$ (the
so-called {\em order conditions})
can be conveniently derived by considering series of linear differential operators.
We denote by $A$ and $B$ the Lie operators
associated with $f^{[a]}$ and $f^{[b]}$, respectively. For each smooth function $g:\R^D \rightarrow \R$, $A\, g$ and $B\, g$ are smooth functions  defined as
\begin{equation*}
  A \, g (x) = \left.\frac{d}{d \tau}\right|_{\tau=0} g(\varphi_{\tau}^{[a]}(x)), \qquad\quad
  B \, g (x) = \left.\frac{d}{d \tau}\right|_{\tau=0} g(\varphi_{\tau}^{[b]}(x)),
\end{equation*}
for each $x \in \mathbb{R}^D$, that is,
\begin{equation}   \label{eq:3b}
   A\, g (x) = f^{[a]}(x) \cdot \nabla g(x),  \qquad
B\, g (x) = f^{[b]}(x) \cdot \nabla g(x).
\end{equation}
The near-integrable Hamiltonian system (\ref{intro.1}) corresponds in this general framework to considering
equation (\ref{eq:1}) with
\[
   x = (q,p), \qquad  f^{[a]}(x) = J \, \nabla H^{[a]}(q,p), \quad \mbox{ and } \quad  f^{[b]}(x) = J \, \nabla H^{[b]}(q,p),
\]
 being $J$ the canonical symplectic matrix. Therefore, for each smooth function $g$ one has
\[
  A\, g= f^{[a]} \cdot \nabla g   =  \sum_j \frac{\partial H^{[a]}}{\partial p_j} \frac{\partial g}{\partial q_j} -  \frac{\partial H^{[a]}}{\partial q_j} \frac{\partial g}{\partial p_j},
\]
and a similar expression for $B\, g$.

It is well known that for any smooth function $g$,
the $\tau$-flow of (\ref{eq:1}) satisfies
\[
   g(\varphi_{\tau}(x)) = \e^{\tau (A +  \eps B)}g (x),
\]
where $\e^{\tau (A +  \eps B)}$ is defined as a series of linear differential operators
\[
   \e^{\tau (A +  \eps B)} = \sum_{k=0}^{\infty} \frac{\tau^k}{k!} (A+\eps B)^k.
\]
The same is true for each part in (\ref{eq:1}):
\begin{equation}\label{eq:4}
  g(\varphi^{[a]}_{\tau}(x))= \e^{\tau \, A}\,  g (x), \qquad \quad
  g(\varphi^{[b]}_{\tau}(x))= \e^{ \tau \, \eps\,  B}\,  g (x).
\end{equation}
Analogously, for the integrator
 $\psi_{\tau}$ in (\ref{eq:3}), one has
\[
 g(\psi_{\tau}(x)) = \Psi(\tau)\,  g (x),
\]
where $\Psi(\tau)$ is a series of linear differential operators defined as
\begin{equation}
 \Psi(\tau) = \e^{a_{1} \tau  A} \, \e^{b_{1}\tau \eps B} \cdots \,
\e^{a_{s}\tau A} \, \e^{b_{s}\tau \eps B} \,
\e^{a_{s+1}\tau A}.    \label{eq:5}
\end{equation}
Notice that the exponentials of  Lie derivatives in (\ref{eq:5}) appear in the reverse order with respect to the maps
in the integrator (\ref{eq:3}).

One of the standard ways of deriving the order conditions for splitting methods is the following. By applying
repeatedly the Baker--Campbell--Hausdorff (BCH) formula \cite{varadarajan84lgl} to the factorization (\ref{eq:5}) corresponding to a consistent (i.e., satisfying (\ref{eq:consistency})) splitting method, one is able to express $\Psi(\tau)$ as the formal exponential of only
one operator:

\begin{equation}
  \label{eq:PsiE}
     \Psi(\tau) = \e^{\tau (A + \eps B + E(\tau,\eps))},
\end{equation}

where
\begin{eqnarray}
\nonumber
E(\tau,\eps) &=&  \tau \, \eps \, p_{ab} [A,B] + \tau^2 \,  \eps \, p_{aba} [[A,B],A] +
\tau^2 \,  \eps^2\, p_{abb} [[A,B],B]) \\
\label{eq:E}
&& + \tau^3 \, \eps\, p_{abaa} [[[A,B],A],A]  + \tau^3 \eps^2 p_{abba} [[[A,B],B],A]\\
&& + \tau^3 \eps^3\, p_{abbb} [[[A,B],B],B]+   \mathcal{O}(\tau^4). \nonumber
\end{eqnarray}
Here the symbol $[A, B]$ stands for the commutator of the Lie operators $A$ and $B$, and
$p_{ab},p_{abb},p_{aba},p_{abbb},\ldots$ are polynomials in the parameters $a_i,b_i$ of the splitting scheme.
In particular,
\[
p_{ab} = \frac12 -\sum_{i=1}^{s} b_i c_i, \qquad p_{aba}=\frac12 \sum_{i=1}^{s} b_i c_i (1-c_i)-\frac{1}{12},
\]
where
\begin{equation}
  \label{eq:ci}
c_i = \sum_{j=1}^{i} a_{j},  \qquad i=1,2,\ldots,s
\end{equation}
and $c_{s+1}=1$. The integrator is of order $r$ if $E(\tau,\eps)$ in (\ref{eq:E}) is of size
$\mathcal{O}(\tau^r)$, so that $\Psi(\tau)$
agrees with
the series of linear operators $\e^{\tau(A + \eps B)}$ of the exact flow up to terms of size $\mathcal{O}(\tau^r)$.
In consequence,
the order conditions read
$p_{ab}=p_{abb}=p_{aba}=\cdots =0$ up
to the order considered. For symmetric methods, $\Psi_{-\tau}=\Psi_{\tau}^{-1}$, and thus $E(-\tau,\eps)=E(\tau,\eps)$, so that
$E(\tau,\eps)$ only involves even powers of $\tau$, that is, $p_{w}=0$ for any word $w$ with an even number of letters in the alphabet $\{a,b\}$.

In (\ref{eq:E}) we have considered the classical Hall basis associated to the Hall words
$a,b,ab,abb,aba,abbb,abba,abaa,\ldots$
\cite{reutenauer93fla}.  The coefficients $p_{w}$ in
(\ref{eq:E}) corresponding to each Hall word $w$ can be
systematically obtained using the results in~\cite{murua06tha} in
terms of rooted trees and iterated integrals.  An efficient algorithm
(based on the results in~\cite{murua06tha}) of the BCH formula and
related calculations that allows one to obtain expression
(\ref{eq:E}) up to terms of arbitrarily high degree is
presented in~\cite{casas08aea}.


\subsection{Generalized order}
\label{ss:gen-order}

We are particularly interested in the manner in which the local error $\psi_{\tau}(x) - \varphi_{\tau}(x)$ decreases as
$\eps \rightarrow 0$. For instance, from the results in the precedent subsection, it is clear that for any consistent symmetric method
the local error satisfies
$\psi_{\tau}(x) = \varphi_{\tau}(x) +
\mathcal{O}(\eps \, \tau^{3})$. Alternatively,
\[
    \Psi(\tau)- \e^{\tau (A + \eps\, B)}=\mathcal{O}(\eps \, \tau^{3}) \quad \mbox{ as } \quad  (\tau,\eps) \rightarrow (0,0).
\]
If in addition $p_{aba} = 0$ in (\ref{eq:E}), then
\[
   \psi_{\tau}(x) = \varphi_{\tau}(x) + \mathcal{O}(\eps \, \tau^{5} + \eps^2 \, \tau^{3}) \quad
   \mbox{ as } \quad  (\tau,\eps) \rightarrow (0,0).
\]
In that case, we say that such a method is of (generalized) order $(4,2)$.
%
More generally, we will say \cite{mclachlan95cmi} that an integration method for the system (\ref{eq:1}) is of generalized
order $(r_1,r_2,\ldots,r_m)$ (where $r_1 \geq r_2 \geq \cdots \geq r_m$) if the local error satisfies that
\begin{equation*}
\psi_{\tau}(x) - \varphi_{\tau}(x) =
\mathcal{O}(\eps \tau^{r_{1}+1}+
  \eps^{2} \tau^{r_{2}+1}+\cdots+\eps^{m}
\tau^{r_{m}+1}).
\end{equation*}
Recall from Subsection~\ref{ss:BCH} that for symmetric integrators, the remainder $E(\tau,\eps)$ in (\ref{eq:PsiE}) is even with respect to $\tau$, and thus the generalized order $(r_1,r_2,\ldots,r_m)$ of symmetric schemes must have even $r_j$.

\subsection{Generalized order conditions}

The conditions that the coefficients $a_i, b_i$ must satisfy for a splitting method to be of a prescribed
(generalized) order $(r_1,r_2,\ldots,r_m)$ can be obtained, of course, by computing the polynomials $p_w$ in expression (\ref{eq:E}) with the BCH formula and then equating each term to zero up to the considered order. Thus, in particular,
a consistent symmetric scheme of order $(6,2)$ requires that $p_{aba} = p_{abaaa} = 0$.

There exist, however, other more systematic procedures to derive these order conditions. In what follows, we present
a strategy that allows us to get in a direct way a set of necessary and sufficient independent order conditions
for generic splitting methods.

As a first step, we consider $Z(\tau) = \e^{\tau(A+\eps \, B)} \e^{-\tau A}$, which is the formal solution of
the initial value problem
\begin{equation}
  \label{eq:odeZ}
  \frac{d}{d \tau} Z(\tau) = \eps\, Z(\tau) \, C(\tau), \qquad Z(0)=I,
\end{equation}
where
\begin{equation}
  \label{eq:C(h)}
  C(\tau) =  \e^{\tau A} B \e^{-\tau A} =  \sum_{n=1}^{\infty} \tau^{n-1} C_n,
\end{equation}
with
\[
 C_1=B, \qquad C_n= \frac{1}{(n-1)!} \, [\underbrace{A,[A,...,[A}_{n-1 \; \mbox{\scriptsize times}},B]]], \quad n > 1.
\]
 On the other hand, applying repeatedly the identity
 $\e^{\tau A} \e^{h \,  B} \e^{-\tau A} = \e^{h \, C(\tau)}$ to eq. (\ref{eq:5}) and taking into account (\ref{eq:ci}),
 we arrive at
\begin{equation}
  \label{eq:8}
  \Psi(\tau) = \widehat Z(\tau)\, \e^{c_{s+1} \tau A}, \qquad \mbox{ where } \qquad
\widehat Z(\tau) = \e^{ \eps\,  b_1 \tau \, C(c_1 \tau)} \cdots \; \e^{ \eps\,  b_s \tau \, C(c_s \tau)}.
\end{equation}
Notice that, if the splitting method is consistent, then $c_{s+1}=1$.
We thus have that a splitting method is of order $(r_1,r_2,\ldots,r_m)$ if and only if $c_{s+1}=1$ and
\begin{equation}
  \label{eq:6}
\widehat Z(\tau)- Z(\tau)=
  \mathcal{O}(\eps \tau^{r_{1}+1}+
  \eps^{2} \tau^{r_{2}+1}+\cdots+\eps^{m}
\tau^{r_{m}+1}).
\end{equation}
We then expand both $Z(\tau)$ and $\tilde{Z}(\tau)$ as power series of $\eps$  and compare their coefficients. First,
applying Neumann iteration to (\ref{eq:odeZ}) we get
\begin{eqnarray*}
  Z(\tau) - I &=&  \eps\,  \int_{0}^{\tau} Z(s_1) \, C(s_1) \, ds_1 \\
        &=&  \eps\, \int_{0}^{\tau} C(s_1) \, ds_1
+ \eps^2\,  \int_{0}^{\tau} \int_{0}^{s_1} C(s_2) \, C(s_1) \, ds_1\, ds_2 \\
&& + \, \eps^3\,  \int_{0}^{\tau} \int_{0}^{s_1} \int_{0}^{s_2} C(s_3) \, C(s_2) \, C(s_1)\, ds_1\, ds_2 \, ds_3 + \cdots \\
&=& \sum_{k\geq 1} \eps^k \sum_{j_1,\ldots,j_k \geq 1} \frac{\tau^{j_1+\cdots +j_k}}{(j_1+\cdots + j_k)\cdots (j_1+j_2) j_1} \, C_{j_1} \cdots C_{j_k},
\end{eqnarray*}
where in the last equality we have introduced explicitly the expression for $C(\tau)$ given by (\ref{eq:C(h)}).
On the other hand, by expanding the exponentials of $\tilde{Z}(\tau)$ in (\ref{eq:8}), we have
\begin{eqnarray*}
\widehat Z(\tau)-I &=& \tau \, \eps \, \sum_{i=1}^{s} b_i C(c_i \tau) \\
&& + \tau^2 \, \eps^2 \left(
\sum_{i=1}^{s} \frac{b_i^2}{2} C(c_i \tau)^2 + \sum_{i=1}^{s-1} \sum_{j=i+1}^{s} b_i b_j C(c_i \tau) C(c_j \tau)
\right) + \cdots \\
&=& \sum_{k\geq 1} \tau^k\eps^k \sum_{1\leq i_1 \leq \cdots \leq i_k \leq s} \frac{b_{i_1}\cdots b_{i_k}}{\sigma_{i_1 \cdots i_k}}\,
C(c_{i_1}\tau) \cdots C(c_{i_k}\tau)\\
&=& \sum_{k\geq 1}\eps^k
\sum_{j_1,\ldots,j_k \geq 1} \tau^{j_1+\cdots +j_k}
\left(\sum_{1\leq i_1 \leq \cdots \leq i_k \leq s} \frac{b_{i_1}\cdots b_{i_k}}{\sigma_{i_1 \cdots i_k}}\,
c_{i_1}^{j_1-1} \cdots c_{i_k}^{j_k-1}\right) \,
 C_{j_1} \cdots C_{j_k},
\end{eqnarray*}
where
\begin{eqnarray*}
  \sigma_{i_1\cdots i_k}=1 &\mbox{ if }& i_1 < \cdots < i_k,\\
  \sigma_{i_1\cdots i_k}=\frac{1}{\ell !}\, \sigma_{i_{\ell+1}\cdots i_k} &\mbox{ if }& i_1=\cdots = i_{\ell} < i_{\ell+1} \leq \cdots \leq i_k.
\end{eqnarray*}
In this way a splitting method is of order $(r_1,\ldots,r_m)$ if and only if
\begin{equation}
\label{eq:ocond}
  \sum_{1\leq i_1 \leq \cdots \leq i_k \leq s} \frac{b_{i_1}\cdots b_{i_k}}{\sigma_{i_1 \cdots i_k}}\,
c_{i_1}^{j_1-1} \cdots c_{i_k}^{j_k-1} = \frac{1}{(j_1+\cdots + j_k)\cdots (j_1+j_2) j_1}
\end{equation}
for each $k=1,\ldots,m$ and each multi-index (i.e., $k$-tuple of positive integers) $(j_1,\ldots,j_k)$
such that $j_1+\cdots+j_k\leq r_k$.

Conditions (\ref{eq:ocond}) (one condition for each multi-index) have been obtained in \cite{thalhammer08hoe}
in the context of order conditions of splitting operators for unbounded operators $A$ and $B$.
Nevertheless, such order conditions are not all independent. For instance, it can be checked that if condition (\ref{eq:ocond}) holds for the multi-indices $(1,2)$, $(2)$, and $(1)$, then the condition for $(2,1)$ is also fulfilled. That kind of dependencies are a consequence of the fact that both $Z(\tau)$ and $\widehat Z(\tau)$ are exponentials of Lie series in the non-commuting indeterminates $C_1,C_2,\ldots$. A set of independent order conditions can be obtained (by virtue of Theorems 3.2 and 6.1 in \cite{reutenauer93fla}) by considering a particular subset of multi-indices,
the so-called {\em Lyndon multi-indices}. Let us
consider the lexicographical order $<$ (i.e., the order used when
ordering words in the dictionary) on the set of multi-indices. A
multi-index $(i_1,\ldots,i_m)$ is a Lyndon multi-index if
$(i_1,\ldots,i_k) < (i_{k+1}, \ldots, i_m)$ for each $1 \leq k <
m$. For instance, the subset of Lyndon multi-indices $(j_1,\ldots,j_k)$ such that $j_1+\cdots+j_k\leq 5$ is
\[\{(1),(2),(3),(4),(5),(1,2),(1,3),(1,4),(2,3),(1,1,2),(1,1,3),(1,2,2),(1,1,1,2)\}.\]

Taking into account these considerations, we finally arrive at the following result.
\begin{theorem}
\label{th:1}
  A splitting method of the form (\ref{eq:3}) is of generalized order $(r_1,\ldots,r_m)$ if and only if $c_{s+1}=1$ and
  (\ref{eq:ocond}) holds for  $k=1,\ldots,m$ and each Lyndon multi-index
  $(j_1,\ldots,j_k)$ such that $j_1+\cdots+j_k\leq r_k$.
For symmetric methods, only Lyndon multi-indices $(j_1,\ldots,j_k)$ with odd $j_1+\cdots+j_k$ need to be considered.
\end{theorem}
For illustration, in Table \ref{tab.1} we collect explicitly conditions (\ref{eq:ocond}) corresponding to some particular multi-indices,
whereas in Table \ref{tab.2} we specify which particular Lyndon multi-indices one has to consider, or equivalently
which conditions (\ref{eq:ocond}) must hold for each
consistent symmetric splitting method of the given generalized order, according to Theorem \ref{th:1}.

\

\begin{table}[h!]
\begin{center}
\begin{tabular}{|l|l|} \hline
  Multi-index  &  \hspace*{0.4cm}  Condition  \\ [0.4ex] \hline
      $(j), \;\;  j \ge 1$  & $\displaystyle \sum_{i=1}^{s} b_i \, c_i^{j-1}=\frac{1}{j}$ \\ [0.7ex]
      $(1,2)$   &  $\displaystyle \sum_{i=1}^{s} \frac{1}{2} b_i^2 c_i + \sum_{1\leq i < j \leq s} b_{i} b_{j} c_{j} = \frac{1}{3}$ \\    [0.6ex]
      $(1,4)$  & $\displaystyle \sum_{i=1}^{s} \frac{1}{2} b_i^2 c_i^3+ \sum_{1\leq i < j \leq s} b_{i} b_{j} c_{j}^3 = \frac{1}{5}$ \\ [0.6ex]
     $(2,3)$ &$ \displaystyle  \sum_{i=1}^{s} \frac{1}{2} b_i^2 c_i^3+ \sum_{1\leq i < j \leq s} b_{i} b_{j} c_{i} c_{j}^2 =
     \frac{1}{10}$ \\ [0.5ex] \hline
\end{tabular}
\end{center}
\caption{\small{Generalized order condition associated with each Lyndon multi-index.}}
\label{tab.1}
\end{table}
\begin{table}[h!]
\begin{center}
\begin{tabular}{|l|l|} \hline
   Generalized order  &  \hspace*{0.4cm} Lyndon multi-indices  \\ \hline
      $(2n,2)$  &  $(3), (5), \ldots, (2n-1)$  \\
      $(8,4)$  &  $(3)$, $(5)$, $(7)$, $(1,2)$ \\
      $(10,4)$  &  $(3)$, $(5)$, $(7)$, $(9)$, $(1,2)$ \\
      $(8,6,4)$  &  $(3)$, $(5)$, $(7)$, $(1,2)$, $(1,4)$, $(2,3)$ \\
      $(10,6,4)$  &  $(3)$, $(5)$, $(7)$, $(9)$, $(1,2)$, $(1,4)$, $(2,3)$ \\ \hline
\end{tabular}
\end{center}
\caption{\small{Lyndon multi-indices corresponding to consistent symmetric splitting methods of a given generalized
order. }}
\label{tab.2}
\end{table}

At this point some remarks must be done. The set of order conditions given by Theorem~\ref{th:1} is
completely equivalent to the order conditions that can be obtained by following the standard approach described in subsection~\ref{ss:BCH}. On the one hand,
the free Lie algebra $\mathcal{L}(C_1,C_2,C_3,\ldots)$ generated by the non-commuting indeterminates $C_1,C_2,C_3,\ldots$,
admits a basis (the Lyndon basis \cite{reutenauer93fla}) in one-to-one correspondence with the set of Lyndon multi-indices. Clearly, if instead of directly comparing the series expansions of $Z(\tau)$ and $\widehat Z(\tau)$ as above, we compare
the formal logarithms $\log(Z(\tau))$ and $\log(\widehat Z(\tau))$, we could obtain one order condition per element in the Lyndon basis. On the other hand, the approach in subsection~\ref{ss:BCH} gives one order condition per element in a basis of
the free Lie algebra $\mathcal{L}(A,B)$ generated by the noncommuting indeterminates $A$ and $B$, and
Lazard elimination theorem~\cite{reutenauer93fla} shows that, as vector spaces, the direct sum of $\mathcal{L}(C_1,C_2,C_3,\ldots)$ with the linear span of $A$ is isomorphic to
$\mathcal{L}(A,B)$ (sending $C_1,C_2,C_3,\ldots$ to $B,[A,B],\frac12 [A,[A,B]],\ldots$ respectively).

\section{New numerical schemes}
\label{sec.3}

There are two different types of symmetric composition schemes (\ref{eq:3}): one in which
the first and last flows correspond to the $A$ part (and thus appropriately called ABA composition),
\begin{equation}   \label{aba.comp}
   \mbox{ ABA: } \quad \varphi^{[a]}_{a_{1} \tau} \circ  \varphi^{[b]}_{b_{1} \tau}\circ \varphi^{[a]}_{a_{2} \tau}\circ
 \cdots \circ \varphi^{[a]}_{a_{2} \tau} \circ  \varphi^{[b]}_{b_{1}\tau} \circ \varphi^{[a]}_{a_{1}\tau}
\end{equation}
and the other in which the role of $\varphi^{[a]}_{\tau}$ and $\varphi^{[b]}_{\tau}$ is interchanged (BAB
composition):
\begin{equation}   \label{bab.comp}
   \mbox{ BAB: } \quad \varphi^{[b]}_{b_{1} \tau} \circ  \varphi^{[a]}_{a_{2} \tau}\circ \varphi^{[b]}_{b_{2} \tau}\circ
 \cdots\circ  \varphi^{[b]}_{b_{2} \tau} \circ \varphi^{[a]}_{a_{2}\tau} \circ \varphi^{[b]}_{b_{1}\tau}.
\end{equation}
Notice that both types of composition are closely related: an $s$-stage BAB method is just an $(s+1)$-stage
ABA scheme with $a_1 = 0$, so that to construct BAB methods one has to solve the same order conditions as for
ABA compositions.
Although $A$ and $B$ are qualitatively different here, and therefore both types of composition may lead in principle
to integrators with different performances, in practice, and for the examples analyzed, we have not found substantial differences, so that in what follows we only consider ABA methods for clarity in the presentation.

Constructing particular methods requires solving polynomial equations (e.g. the order conditions of Table \ref{tab.1}
for methods of Table \ref{tab.2}), a problem
whose complexity grows enormously with the number of equations and variables involved. This task can be handled
by computer algebra systems when this number is relatively low. In that case one is
able to get all the solutions and select the one that verifies some previously fixed optimization criterion, such as
minimizing error terms at higher orders and the sum of the absolute value of  the coefficients. In practice, we have
followed this procedure when there are no free parameters and the number of equations to be solved is at most seven.
When this number is larger than seven or there are additional parameters, another strategy based on homotopy
continuation methods has been applied. In the appendix we provide a detailed treatment of the
procedure for a particular method.

\subsection{New methods in the ABA class}

Symmetric schemes of generalized order $(2n,2)$ can be obtained just by solving, in addition to
consistency,  the order conditions
corresponding to the Lyndon multi-indices $(3), (5), \ldots, (2n-1)$ (first line in Table \ref{tab.1}). These equations
result from approximating the integral $\int_0^{\tau} C(s) ds$ in the expression of $Z(\tau)$ by the
quadrature rule
\[
  \sum_{i=1}^s b_i \, C(c_i \tau) = \sum_{i=1}^s b_i \sum_{j \ge 1} c_i^{j-1} \tau^{j-1} C_j
\]
in the expansion of $\tilde{Z}(\tau)$.
Equivalent order conditions were previously derived in \cite{laskar01hos,mclachlan95cmi,mclachlan96mos}, and so the
same methods are obtained here.  Methods in this family have all their coefficient positive and good
stability properties.
In the tests carried out in this paper we will take the most efficient
ABA scheme of order $(8,2)$ for comparison, which we denote by \texttt{ABA82}.

\paragraph{Generalized order (10,4).} According with Table \ref{tab.2}, there are five order conditions in addition to
consistency (\ref{eq:consistency}), for a total number of seven equations to be satisfied by the coefficients. As a 
consequence, the minimum number of stages is six. A more efficient method can be obtained, however, by taking
an additional stage and choosing the corresponding free parameter to reduce the error terms at a higher order.
The sequence of coefficients is
\begin{equation}   \label{10.4}
    a_1 \, b_1 \, a_2 \, b_2 \,  a_3 \, b_3 \, a_4 \, b_4 \, a_4 \, b_3 \,  a_3 \, b_2 \, a_2 \, b_1 \, a_1
\end{equation}
and their values are collected in Table~\ref{Tabla3b} (method denoted
by \texttt{ABA104}). Observe that, as expected, one of the $a_i$  and
one of the $b_j$ coefficients are negative (it is known that this feature is unavoidable
for any splitting
method of order higher than two~\cite{goldman96noo,sheng89slp,suzuki91gto}), but they have a relatively
small absolute value.

\paragraph{Generalized order (8,6,4).} Here we have, in addition to consistency, six order conditions, for a total
of eight equations, so that the minimum number of stages is seven. Hence the sequence of coefficients for the resulting methods is also 
as in (\ref{10.4}). There are 30 real solutions, and the
one referred to as method \texttt{ABA864} in Table~\ref{Tabla3b} minimizes the sum of the absolute values of its
coefficients.

\paragraph{Generalized order (10,6,4).} An additional stage is required in this case to verify the order condition
associated with multi-index (9) in Table \ref{tab.1}. Therefore, the minimum number of stages is eight,
with sequence
\begin{equation}   \label{10.6.4}
    a_1 \, b_1 \, a_2 \, b_2 \,  a_3 \, b_3 \, a_4 \, b_4 \, a_5 \, b_4 \, a_4 \, b_3 \,  a_3 \, b_2 \, a_2 \, b_1 \, a_1.
\end{equation}
By following the construction strategy exposed in the appendix, we have obtained several solutions for the
order conditions. Among those possessing reasonably small coefficients, we have selected the solution with the
smallest leading terms of the local error. This corresponds to method \texttt{ABA1064} in Table~\ref{Tabla3b}.

\begin{table}[htb]
\centering
\begin{tabular}{cccl}
\hline\noalign{\smallskip}
{\tt id} & {\tt order} & {\tt stages}& \hspace*{3.5cm} $a_i,  b_i$ \\
\noalign{\smallskip}\hline\noalign{\smallskip}
{\small {\tt ABA104}} & {\small $(10,4)$} & {\small 7} &
{\footnotesize \begin{tabular}{rcl}
$a_1$ &=& {\tt \ 0.04706710064597250612947887637243678556564} \\
$a_2$ &=& {\tt \ 0.1847569354170881069247376193702560968574} \\
$a_3$ &=& {\tt \ 0.2827060056798362053243616565541452479160} \\
$a_4$ &=& {\tt  -0.01453004174289681837857815229683813033908} \\
$b_1$ &=& {\tt \ 0.1188819173681970199453503950853885936957}\\
$b_2$ &=& {\tt \ 0.2410504605515015657441667865901651105675}\\
$b_3$ &=& {\tt  -0.2732866667053238060543113981664559460630}\\
$b_4$ &=& {\tt \ 0.8267085775712504407295884329818044835997}\\
    \end{tabular} }\\
\noalign{\smallskip}\hline\noalign{\smallskip}
{\small {\tt ABA864}} & {\small $(8,6,4)$} & {\small 7} &
{\footnotesize \begin{tabular}{rcl}
$a_1$ &=& {\tt \ 0.0711334264982231177779387300061549964174} \\
$a_2$ &=& {\tt \ 0.241153427956640098736487795326289649618} \\
$a_3$ &=& {\tt \ 0.521411761772814789212136078067994229991} \\
$a_4$ &=& {\tt  -0.333698616227678005726562603400438876027} \\
$b_1$ &=& {\tt \ 0.183083687472197221961703757166430291072} \\
$b_2$ &=& {\tt \ 0.310782859898574869507522291054262796375} \\
$b_3$ &=& {\tt  -0.0265646185119588006972121379164987592663} \\
$b_4$ &=& {\tt \ 0.0653961422823734184559721793911134363710} \\
   \end{tabular} }\\
\noalign{\smallskip}\hline\noalign{\smallskip}
{\small {\tt ABA1064}} & {\small $(10,6,4)$} & {\small 8} &
{\footnotesize \begin{tabular}{rcl}
$a_1$ &=& {\tt \ 0.03809449742241219545697532230863756534060} \\
$a_2$ &=& {\tt \ 0.1452987161169137492940200726606637497442} \\
$a_3$ &=& {\tt \ 0.2076276957255412507162056113249882065158} \\
$a_4$ &=& {\tt \ 0.4359097036515261592231548624010651844006} \\
$a_5$ &=& {\tt  -0.6538612258327867093807117373907094120024} \\
$b_1$ &=& {\tt \ 0.09585888083707521061077150377145884776921} \\
$b_2$ &=& {\tt \ 0.2044461531429987806805077839164344779763} \\
$b_3$ &=& {\tt \ 0.2170703479789911017143385924306336714532} \\
$b_4$ &=& {\tt  -0.01737538195906509300561788011852699719871} \\
  \end{tabular} }\\
\noalign{\smallskip}\hline
\end{tabular}
\caption{\small{
Coefficients for ABA symmetric splitting methods of
generalized order $(10,4)$, $(8,6,4)$ and $(10,6,4)$.}}
\label{Tabla3b}
\end{table}

\subsection{A simple example}

To illustrate the efficiency of these schemes we take the perturbed Kepler problem with
Hamiltonian
\begin{equation}
H=\frac{\displaystyle1}{\displaystyle2}(p_{1}^{2}+p_{2}^{2})-\frac
{\displaystyle1}{\displaystyle r}-\frac{\displaystyle\eps
}{\displaystyle2r^{3}}\,\left(  1-\frac{\displaystyle3q_1^{2}%
}{\displaystyle r^{2}}\right),  \label{ne.2}%
\end{equation}
where $r=\sqrt{q_1^{2}+q_2^{2}}$. This Hamiltonian is a first approximation used to describe
the dynamics of a satellite moving into the
gravitational field produced by a slightly oblate spherical
planet and whose motion takes place in a plane containing the
symmetry axis of the planet \cite{meirovich88moa}.


We consider this simple problem to test the relative performance of the
methods obtained in this work in comparison with schemes presented in
\cite{laskar01hos,mclachlan95cmi}. The following schemes are used:
\begin{itemize}
\item \texttt{ABA82}: The 4-stage (8,2) ABA method given in \cite{laskar01hos,mclachlan95cmi}.
\item \texttt{ABA84}: The 5-stage (8,4) ABA method of \cite{mclachlan95cmi}.
\item \texttt{ABA104}: The 7-stage (10,4) method given in Table \ref{Tabla3b}.
\item \texttt{ABA864}: The 7-stage (8,6,4) method of Table \ref{Tabla3b}.
\item \texttt{ABA1064}: The 8-stage (10,6,4) method whose coefficients are collected in Table \ref{Tabla3b}.
\end{itemize}

For the numerical experiments we take as initial conditions
$q_1=1-e$, $q_2=0$, $p_{1}=0$, $p_{1}=\sqrt{(1+e)/(1-e)}$, with
$e=1/4$, which would correspond to the eccentricity of the
unperturbed Kepler problem. For this system, the strength of the
perturbation depends both on the choice of the small parameter,
$\eps$, and the initial conditions. We
integrate along the interval $t\in[0,10000]$ and compute the
averaged error in energy as well as the averaged error in position and momenta
(measured in the 2-norm) of the numerical solutions
evaluated at $t_{k}=20\cdot k, \ k=1,2,\ldots,500$. We take as the
exact solution an accurate approximation obtained using a high
order method with a sufficiently small time step.
%
This numerical test is repeated several times for each method
using different time steps (changing the computational cost for
the numerical integration). Finally, we plot the average errors
versus the time step scaled by the number of stages per step,
i.e. $\tau/s$, in double logarithmic scale, to show how the error
depends on the computational cost (the cost is inversely proportional
to $\tau/s$, and the best methods should provide a given accuracy
with the largest value of $\tau/s$).

\begin{figure}[h!]
\centering
\includegraphics[width = 0.49\textwidth]{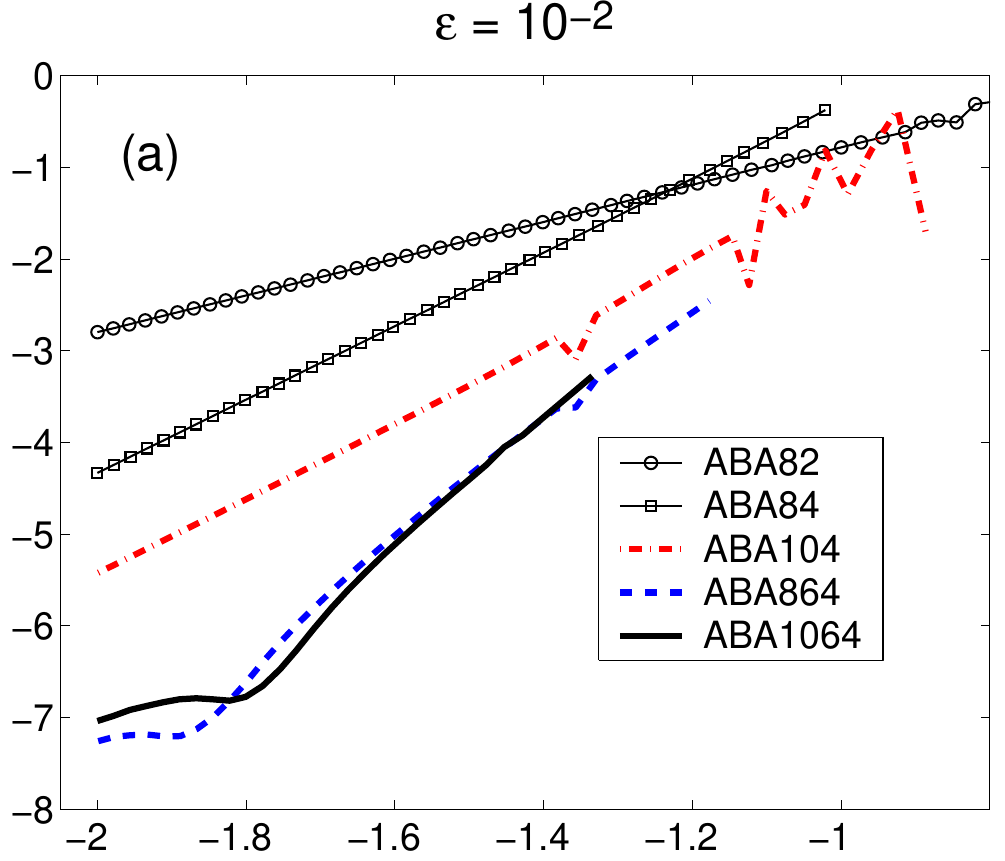}
\includegraphics[width = 0.49\textwidth]{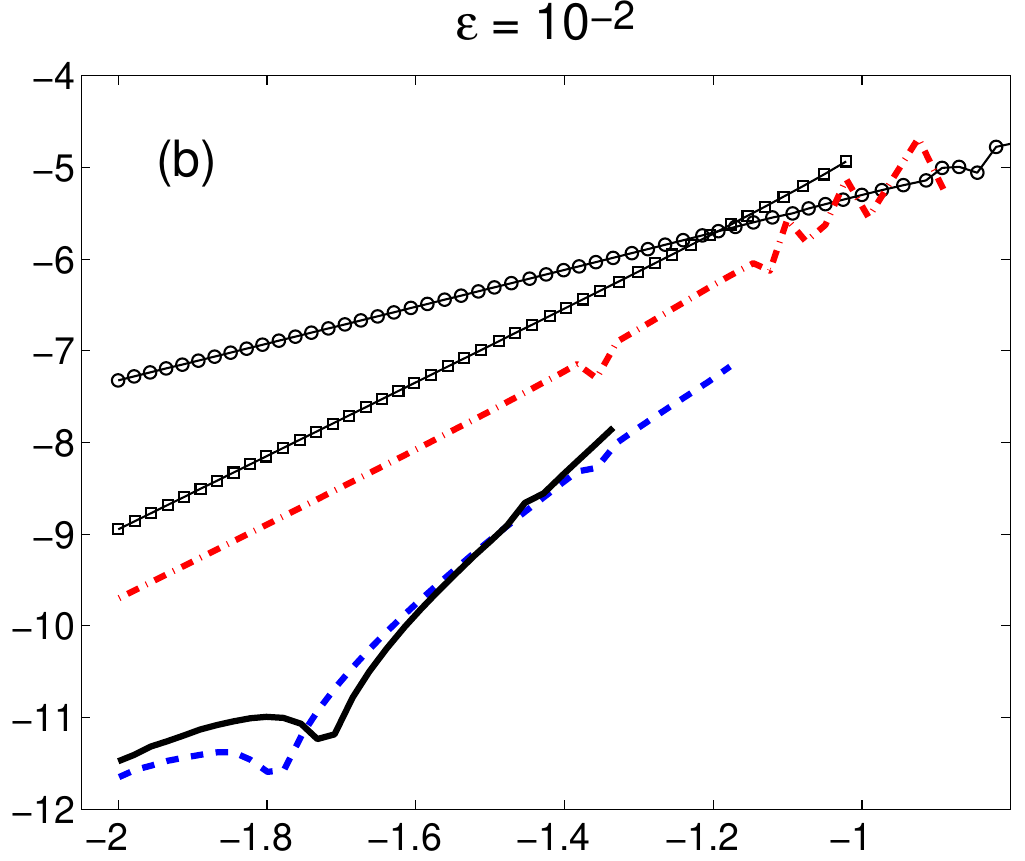}
\includegraphics[width = 0.49\textwidth]{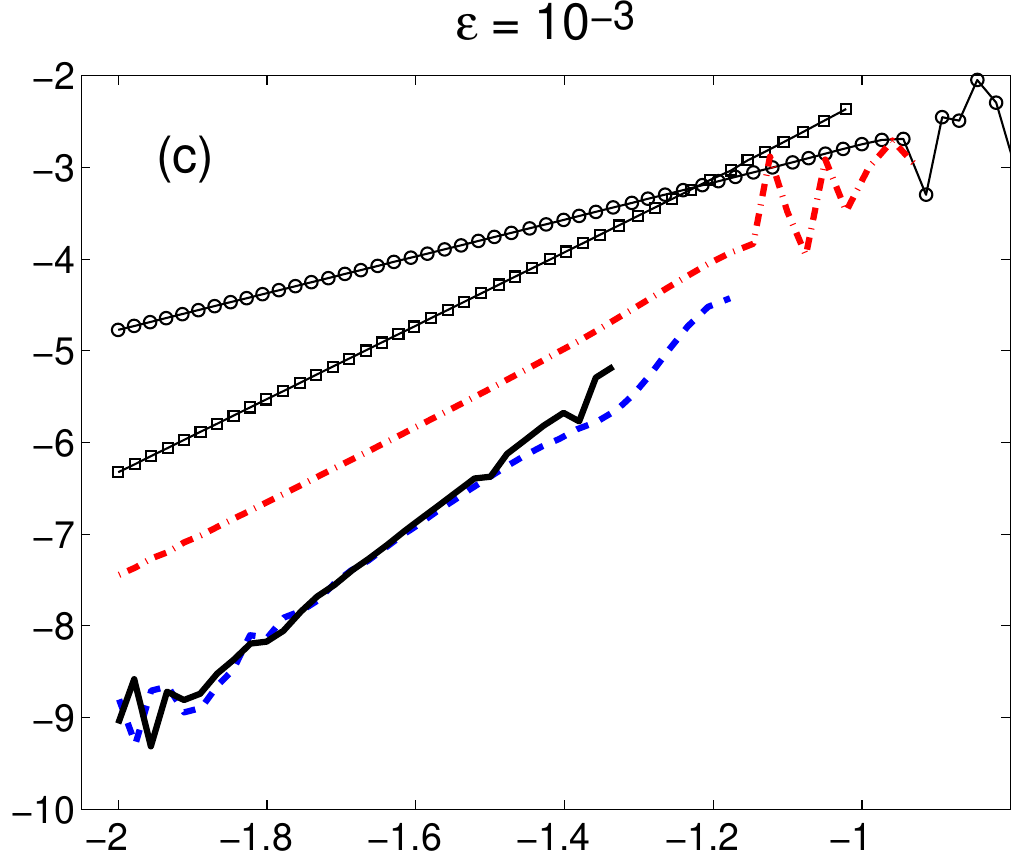}
\includegraphics[width = 0.49\textwidth]{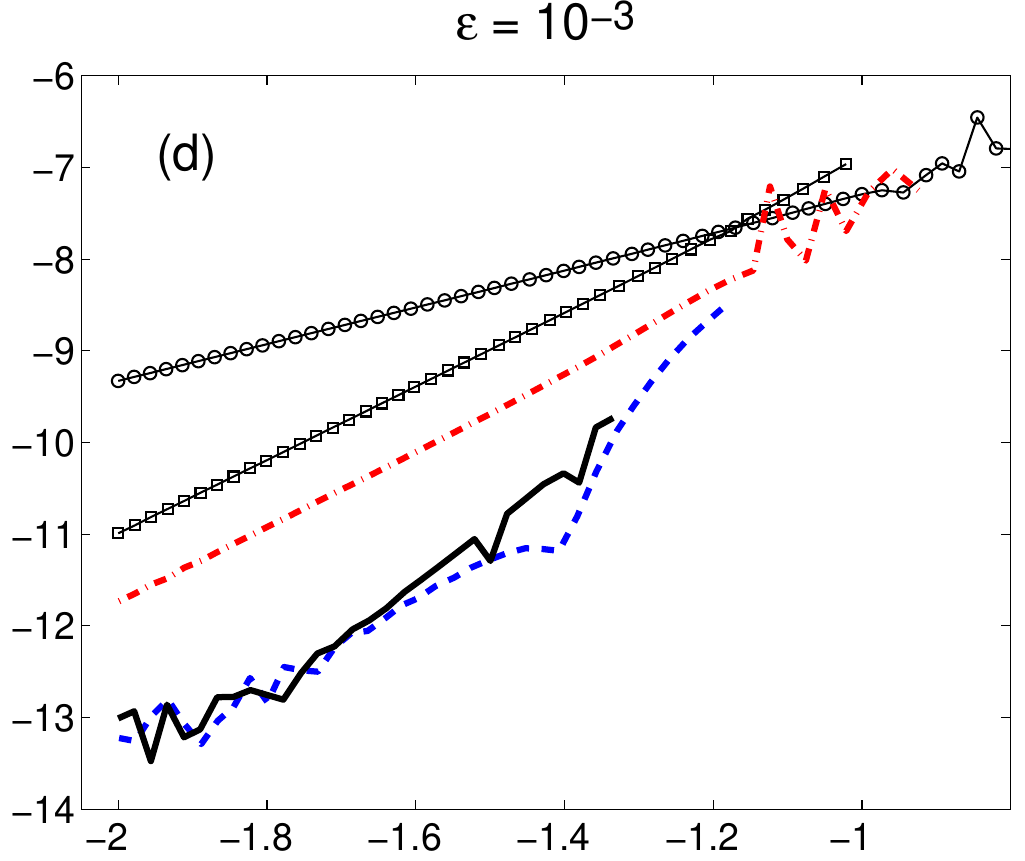}
\caption{\small{Average error in positions and momenta (panels (a) and (c))
and average error in energy (panels (b) and (d)) versus the scaled time
step, $\tau/s$, in a double logarithmic scale for the numerical
integration of the Hamiltonian system (\ref{ne.2}) along the time
interval $t\in[0,10000]$ and measured at times $t_{k}=20\cdot k, \
k=1,2,\ldots,500$.}}
 \label{figEj2-a}
\end{figure}

Figure~\ref{figEj2-a} shows the results obtained for
$\eps=10^{-2},10^{-3}$. In diagrams (a) and (c)  we show the
average error in positions and momenta, whereas in pictures (b) and (d)
we measure the
average error in energy. Notice how the new methods collected in
Table \ref{Tabla3b} are clearly more efficient than \texttt{ABA82} and
\texttt{ABA84}.

It should be stressed that, although the coefficients in Table \ref{Tabla3b} have 40 digits of accuracy,
the results displayed in Figure~\ref{figEj2-a} have been obtained, for the sake of illustration,
with a standard Fortran compiler in double precision. The code generating the results for the averaged error
in energy is available at the website \texttt{www.gicas.uji.es/software.html}.

\subsection{Splitting methods with approximate flows}

 We have so far assumed that the exact $\tau$-flow maps $\varphi_{\tau}^{[a]}$  and $\varphi_{\tau}^{[b]}$ in the splitting method (\ref{eq:3}) are both available. This is the case for the simple example (\ref{ne.2}) considered in the precedent subsection, where the Hamiltonian is split as the sum of a Keplerian Hamiltonian $H^{[a]}(q,p)$ and a perturbation $\eps H^{[b]}(q)$ that only depends on the positions.
However, if instead of the perturbation of example (\ref{ne.2}), one has a perturbation that depends on both positions $q$ and momenta $p$, then,
in general, 
the exact $\varphi_{\tau}^{[b]}$ will no longer be available.
In that case, instead of the splitting method (\ref{eq:3}), we will consider a composition of the form
\begin{equation}   \label{eq:3approx}
\widetilde{\psi}_{h} = \varphi^{[a]}_{a_{s+1} \tau} \circ  \widetilde{\varphi}^{[b]}_{b_{s} \tau}\circ \varphi^{[a]}_{a_{s} \tau}\circ
 \cdots \circ \varphi^{[a]}_{a_{2} \tau} \circ  \widetilde{\varphi}^{[b]}_{b_{1}\tau} \circ \varphi^{[a]}_{a_{1}\tau},
\end{equation}
where $\widetilde{\varphi}^{[b]}_{\tau}$ is an approximation of $\varphi^{[b]}_{\tau}$ obtained by applying some numerical integrator to the Hamiltonian $\eps H^{[b]}(q,p)$.

In what follows, we assume that $\widetilde{\varphi}^{[b]}_{\tau}$ represents one step of some 2nd-order symmetric method.  In that case, the series of differential operators corresponding to
$\widetilde{\varphi}^{[b]}_{\tau}$ is of the form
$$
\widetilde{\Phi}^{[b]}_{\tau}=\e^{\tau \eps B + (\tau \eps)^{3} D_3 + (\tau \eps)^{5} D_5 +\cdots}
$$
instead of just $\e^{\tau \eps B}$, and thus,
the series $\widetilde{\Psi}_h$  of differential operators corresponding to the method (\ref{eq:3approx}) can be obtained from (\ref{eq:5}) by replacing each $\e^{b_j \tau \eps B}$ by
$$
\widetilde{\Phi}^{[b]}_{b_j \tau} =\e^{b_{j}\tau \eps B + (b_{j}\tau \eps)^{3} D_3 + (b_{j}\tau \eps)^{5} D_5 +\cdots}.
$$
%
In consequence, the leading term of the difference $\widetilde{\Psi}_{\tau}-\Psi_{\tau}$ of the respective series corresponding to methods
(\ref{eq:3approx}) and (\ref{eq:3}) is
$$
\left(\sum_{j=1}^{s} b_j^3\right) \eps^3 \tau^3 D_3.
$$
It is then natural to impose, in addition to the generalized order conditions obtained in subsection~\ref{ss:gen-order}, the condition
\begin{equation}   \label{helio.5}
    \sum_{i=1}^s b_i^3 = 0,
\end{equation}
with the aim of reducing the effect of replacing $\varphi^{[b]}_{\tau}$ by $\widetilde{\varphi}^{[b]}_{\tau}$ in (\ref{eq:3}).

The order conditions of scheme (\ref{eq:3approx}) ($\widetilde{\varphi}^{[b]}_{\tau}$
being one step of arbitrary second order symmetric integrator applied to $y'=\eps f^{[b]}(y)$)
can be systematically obtained by generalizing the approach presented in section 2, just by replacing $\widehat Z(\tau)$ in (\ref{eq:8}) by
\[
\widetilde Z(\tau) = \e^{ G(\eps\,  b_1 \tau, c_1 \tau)} \cdots \; \e^{ G(\eps\,  b_s \tau, c_s \tau)},
\]
where
\begin{eqnarray*}
  G(\sigma,\tau) &=&  \e^{\tau A} (\sigma\, B + \sigma^{3} D_3 + \sigma^{5} D_5 +\cdots) \e^{-\tau A} \\
&=&  \sum_{n=1}^{\infty} \tau^{n-1} (\sigma\, C_n + \sigma^{3} E_{3,n} + \sigma^{5} E_{5,n} +\cdots),
\end{eqnarray*}
with
\[
 E_{i,1}=D_i, \qquad E_{i,n}= \frac{1}{(n-1)!} \, [\underbrace{A,[A,...,[A}_{n-1 \; \mbox{\scriptsize times}},D_{i}]]], \quad n > 1.
\]
Thus, the scheme (\ref{eq:3approx}) will have generalized order $(r_1,r_2,r_3,\ldots,r_m)$ if (\ref{eq:3}) is of generalized order $(r_1,r_2,r_3,\ldots,r_m)$ (see subsection~\ref{ss:gen-order}) and in addition,
 $$
\widetilde{Z}(\tau) -\widehat{Z}(\tau) =  \mathcal{O}(\eps^3 \tau^{r_{3}+1}+
  \eps^{4} \tau^{r_{4}+1}+\cdots+\eps^{m}\tau^{r_{m}+1}).
$$
 In particular, we have that
\begin{eqnarray*}
  \widetilde{Z}(\tau) -\widehat{Z}(\tau) &=&  
\left(\sum_{j=1}^{s} b_j^3\right) \eps^3 \tau^3 D_3 + 
\left(\sum_{j=1}^{s} b_j^3c_j\right) \eps^3 \tau^4 [A,D_3] \\
&& + 
\frac12 \left(\sum_{j=1}^{s} b_j^3c_j^2\right) \eps^3 \tau^5 [A,[A,D_3]] + \mathcal{O}(\eps^3 \tau^6 + \eps^4 \tau^4). 
\end{eqnarray*}
This shows that, if the method (\ref{eq:3approx}) is applied with the coefficients of a standard symmetric splitting method (\ref{eq:3}) of generalized order $(r,4)$, the resulting method has generalized order $(r,4,2)$. If the additional condition (\ref{helio.5}) holds, then (\ref{eq:3approx}) recovers the generalized order $(r,4)$ (recall that generalized order  $(r_1,r_2,\ldots,r_m)$ of symmetric methods have even $r_j$).

We have constructed several symmetric methods
of ABA-type
$$
  \varphi^{[a]}_{a_{1} \tau} \circ  \widetilde{\varphi}^{[b]}_{b_{1} \tau}\circ \varphi^{[a]}_{a_{2} \tau}\circ
 \cdots \circ \varphi^{[a]}_{a_{2} \tau} \circ  \widetilde{\varphi}^{[b]}_{b_{1}\tau} \circ \varphi^{[a]}_{a_{1}\tau}
$$
(with $a_1 \neq 0$) satisfying the additional condition (\ref{helio.5}).
Notice that ABA-type compositions are more
convenient than BAB-type methods, since the last stage of the method in the current step can be concatenated with the first stage at the next
step. This is not possible with BAB compositions (\ref{bab.comp}), because
\[
   \widetilde{\varphi}_{b_i \tau}^{[b]}  \circ \widetilde{\varphi}_{b_j \tau}^{[b]} \ne \widetilde{\varphi}_{(b_i + b_j) \tau}^{[b]}.
\]

By following a similar strategy as for the methods collected in section \ref{sec.3}, we have constructed
symmetric schemes within this family of generalized order (8,4), (8,6,4) and (10,6,4), all of them involving the
minimum number of stages. The corresponding coefficients are collected in Table \ref{TAB_ABAH_high}.
For schemes (8,4) we have found all the real solutions and selected the solution that
minimizes the sum of the absolute values of the coefficients (method \texttt{ABAH844}). This method has the
following structure
\[
    \psi_{\tau} = \varphi^{[a]}_{a_{1} \tau} \circ  \widetilde{\varphi}^{[b]}_{b_{1} \tau}\circ \varphi^{[a]}_{a_{2} \tau}\circ
      \widetilde{\varphi}^{[b]}_{b_{2} \tau} \circ \varphi^{[a]}_{a_{3} \tau} \circ  \widetilde{\varphi}^{[b]}_{b_{3} \tau} \circ
      \varphi^{[a]}_{a_{4} \tau} \circ \widetilde{\varphi}^{[b]}_{b_{3} \tau} \circ \varphi^{[a]}_{a_{3} \tau} \circ \widetilde{\varphi}^{[b]}_{b_{2} \tau}
      \circ \varphi^{[a]}_{a_{2} \tau} \circ \widetilde{\varphi}^{[b]}_{b_{1} \tau} \circ \varphi^{[a]}_{a_{1} \tau}.
\]
The procedure for
constructing method \texttt{ABAH1064} is detailed in the appendix, and a similar strategy has been used to
build scheme {\tt ABAH864}.

\begin{table}[ht]
\centering
\begin{tabular}{cccl}
\hline\noalign{\smallskip}
{\tt id} & {\tt order} & {\tt stages} & \hspace*{2cm} $a_i, b_i$ \\
\noalign{\smallskip}\hline\noalign{\smallskip}
{\tt ABAH844} & $(8,4)$ & 6 &
{\footnotesize \begin{tabular}{rcl}
$a_1$ &=& {\tt \ 0.2741402689434018761640565440378637101205} \\
$a_2$ &=& {\tt  -0.1075684384401642306251105297063236526845} \\
$a_3$ &=& {\tt -0.04801850259060169269119541715084750653701} \\
$a_4$ &=& {\tt \ 0.7628933441747280943044988056386148982021} \\
$b_1$ &=& {\tt \ 0.6408857951625127177322491164716010349386} \\
$b_2$ &=& {\tt  -0.8585754489567828565881283246356000103664} \\
$b_3$ &=& {\tt \ 0.7176896537942701388558792081639989754277} \\
\end{tabular} } \\
\noalign{\smallskip}\hline\noalign{\smallskip}
{\tt ABAH864} & $(8,6,4)$ & 8 &
{\footnotesize \begin{tabular}{rcl}
$a_1$ &=& {\tt \ 0.06810235651658372084723976682061164571212} \\
$a_2$ &=& {\tt \ 0.2511360387221033233072829580455350680082} \\
$a_3$ &=& {\tt -0.07507264957216562516006821767601620052338} \\
$a_4$ &=& {\tt  -0.009544719701745007811488218957217113269121} \\
$a_5$ &=& {\tt \ 0.5307579480704471776340674235341732001443} \\
$b_1$ &=& {\tt \ 0.1684432593618954534310382697756917558148} \\
$b_2$ &=& {\tt \ 0.4243177173742677224300351657407231801453} \\
$b_3$ &=& {\tt  -0.5858109694681756812309015355404036521923} \\
$b_4$ &=& {\tt \ 0.4930499927320125053698281000239887162321} \\
\end{tabular} }\\
\noalign{\smallskip}\hline\noalign{\smallskip}
{\tt ABAH1064} & $(10,6,4)$ & 9 &
{\footnotesize \begin{tabular}{rcl}
$a_1$ &=& {\tt \ 0.04731908697653382270404371796320813250988} \\
$a_2$ &=& {\tt \ 0.2651105235748785159539480036185693201078} \\
$a_3$ &=& {\tt  -0.009976522883811240843267468164812380613143} \\
$a_4$ &=& {\tt  -0.05992919973494155126395247987729676004016} \\
$a_5$ &=& {\tt \ 0.2574761120673404534492282264603316880356} \\
$b_1$ &=& {\tt \ 0.1196884624585322035312864297489892143852} \\
$b_2$ &=& {\tt \ 0.3752955855379374250420128537687503199451} \\
$b_3$ &=& {\tt  -0.4684593418325993783650820409805381740605} \\
$b_4$ &=& {\tt \ 0.3351397342755897010393098942949569049275} \\
$b_5$ &=& {\tt \ 0.2766711191210800975049457263356834696055} \\
   \end{tabular} }\\
\noalign{\smallskip}\hline
\end{tabular}
\caption{\small{
Coefficients for ABA symmetric splitting methods of generalized order (8,4), (8,6,4) and (10,6,4)
especially adapted to be used when the flow $\varphi^{[b]}_{\tau}$ is approximated by a
symmetric 2nd-order method $\widetilde{\varphi}^{[b]}_{\tau}$. This happens, in particular, 
when Heliocentric coordinates are used for the integration of the Solar System, as it is shown in
section \ref{sec.4}.
}}
\label{TAB_ABAH_high}
\end{table}

\section{Application to the integration of the Solar System}
\label{sec.4}

In this section we illustrate how the new families of methods proposed here 
behave when they are used in the numerical integration of the simplest model of the Solar System, i.e.,
a main massive body (the Sun) and a set of particles (the planets) orbiting the Sun
following almost Keplerian trajectories. It is not our intention to carry out a detailed
treatment of this problem, but rather to check the performance of the new methods
and compare them with other well established
schemes designed for near-integrable Hamiltonian systems such as those presented in 
\cite{laskar01hos} and \cite{mclachlan95cmi}. We instead refer the reader to reference 
\cite{farres12sif}, where this issue is handled in much more
detail.

As stated in the Introduction, integrating numerically the gravitational
N-body problem requires first to choose a convenient set of
canonical coordinates. Two widely used coordinate systems where the corresponding
Hamiltonian (\ref{n-body}) adopts the form (\ref{intro.1}), suitable 
to the application of the integration schemes developed in this work, are
Jacobi and Heliocentric coordinates. In the former, the position of each
planet is taken relative to the barycenter of the previous $i$ bodies, whereas
in the later the position of each planet is taken with respect to the Sun.

Methods of Table \ref{Tabla3b} are particularly appropriate for long time integrations of the N-body
problem in Jacobi coordinates, and extensive numerical experiments with different planetary
configurations have been carried out in \cite{farres12sif}. Here we will restrict ourselves to Heliocentric
coordinates.

In this set the coordinates $\mathbf{r}_i$ are the relative positions of each planet with respect to the Sun:
\begin{equation}  \label{helio.1}
   \mathbf{r}_0 = \mathbf{q}_0, \qquad \mathbf{r}_i = \mathbf{q}_i - \mathbf{q}_0, \qquad i=1,\ldots, n
\end{equation}
whereas  the conjugate momenta read 
\begin{equation}  \label{helio.2}
   \mathbf{\tilde{r}}_0  = \mathbf{p}_0 + \cdots + \mathbf{p}_n, \qquad\quad  \mathbf{\tilde{r}}_i = \mathbf{p}_i
\end{equation}
and the Hamiltonian (\ref{n-body}) is given by \cite{laskar90lmm}
\begin{equation}   \label{helio.3}
  H_{\mbox{\scriptsize He}}  =  \sum_{i=1}^n    \left( \frac{1}{2} \| \mathbf{\tilde{r}}_i\|^2 \, \frac{m_0 + m_i}{m_0 m_i}
      - G \frac{m_0 m_i}{\| \mathbf{r}_i\|} \right) + 
        \sum_{0 < i < j \le n} \left( \frac{\mathbf{\tilde{r}}_i \cdot  \mathbf{\tilde{r}}_j}{m_0} - G \frac{m_i m_j}{\Delta_{ij}} \right),  
\end{equation}        
where $\Delta_{ij} = \| \mathbf{r}_i - \mathbf{r}_j \|$ for $i,j > 0$.

Heliocentric coordinates have 
the advantage, compared with Jacobi coordinates, that adding a new body to the model does not change
the origin and the new Hamiltonian is easily updated. On the other hand, the perturbation $H^{[b]}$ depends on
both positions and momenta and is not integrable by itself. Hence, when considering the splitting (\ref{eq:1}), the
equation $x' = f^{[b]}(x)$ is not exactly solvable. Nevertheless, if the corresponding flow $\varphi_{\tau}^{[b]}$ is approximated by
a 2nd-order method, we can use the new splitting methods of Table \ref{TAB_ABAH_high}. Notice that
$H^{[b]}$ is the sum of two terms, one depending only on positions and the other
depending only on momenta
\[
    H^{[b]}(\mathbf{r}, \mathbf{\tilde{r}}) =   H^{[b_{\mathrm{a}}]}(\mathbf{\tilde{r}}) + H^{[b_{\mathrm{b}}]}(\mathbf{r}) = 
     \sum_{0 < i < j \le n} \left( \frac{\mathbf{\tilde{r}}_i \cdot  \mathbf{\tilde{r}}_j}{m_0} - G \frac{m_i m_j}{\Delta_{ij}} \right),  
\]
so each equation $x' = f^{[b_i]}(x)$, with $f^{[b_i]}(x) = J \nabla H^{[b_i]}(x)$, is exactly solvable with
flow $\varphi_{\tau}^{[b_i]}$. Therefore, we can approximate $\varphi_{\tau}^{[b]}$ by 
the second order symmetric scheme
\begin{equation}  \label{helio.4}
  \widetilde{\varphi}_{b_i \tau}^{[b]} = \varphi^{[b_{\mathrm{a}}]}_{b_i \tau/2}\circ\varphi^{[b_{\mathrm{b}}]}_{b_i \tau}\circ
     \varphi^{[b_{\mathrm{a}}]}_{b_i \tau/2}.
\end{equation}
Another possibility consists in taking splitting methods of the form
\[
    \psi_{\tau} = \varphi^{[a]}_{a_{s+1} \tau} \circ  \varphi^{[b_{\mathrm{b}}]}_{c_{s}\tau} \circ \varphi^{[b_{\mathrm{a}}]}_{b_{s}\tau}\circ 
    \varphi^{[a]}_{a_{s} \tau}\circ
 \cdots\circ  \varphi^{[b_{\mathrm{b}}]}_{c_{1}\tau} \circ \varphi^{[b_{\mathrm{a}}]}_{b_{1}\tau} \circ \varphi^{[a]}_{a_{1}\tau}
\]
and obtaining the appropriate coefficients $a_i$, $b_i$, $c_i$ satisfying the required order conditions. These 
can be derived, for instance, by analyzing the free Lie algebra generated by the three Lie derivatives corresponding
to each piece of the Hamiltonian. The number of order conditions grows rapidly with the order, however, 
in comparison with
splitting schemes involving only two parts. In particular, 
(10,6,4) methods require, in addition to consistency, 22 order conditions, while methods of the form
(\ref{eq:3approx}) with (\ref{helio.4}) only need to satisfy 8 order conditions. This being the case, in what follows we consider
splitting methods of the form (\ref{eq:3approx}), with the approximate flow $\widetilde{\varphi}_{\tau}^{[b]}$ given by the leapfrog composition
(\ref{helio.4}).

\begin{figure}[h!]
\centering
\includegraphics[width = 0.49\textwidth]{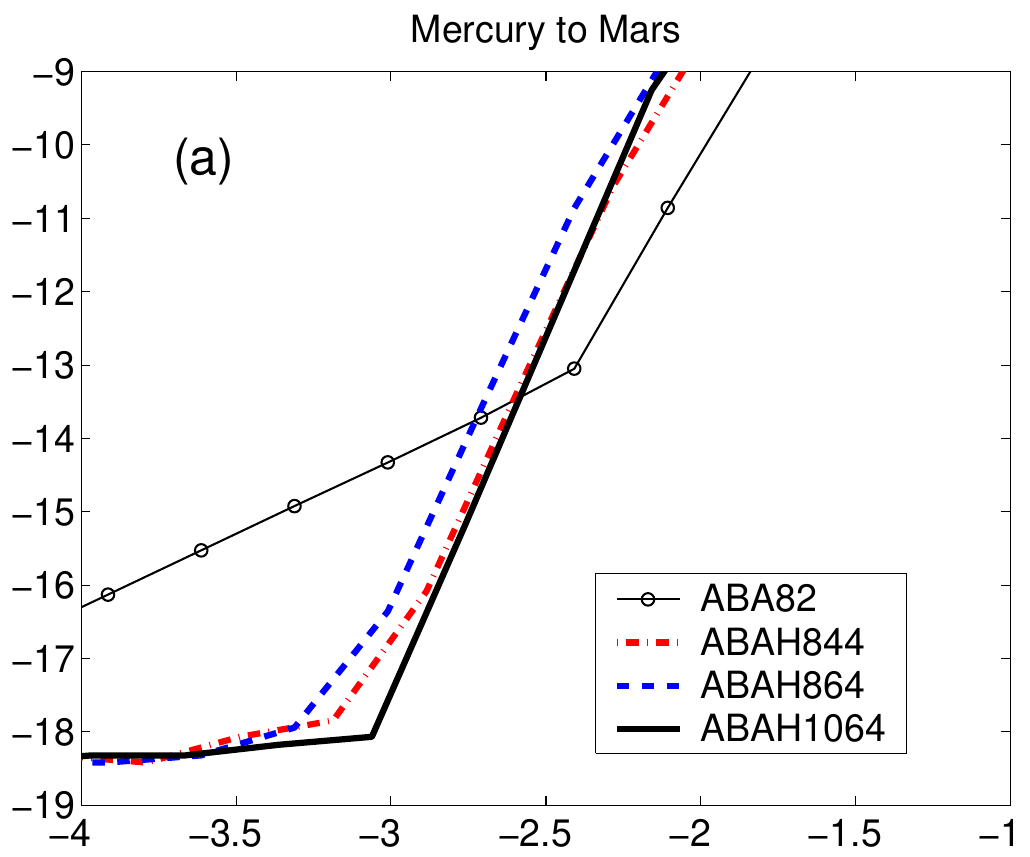}
\includegraphics[width = 0.49\textwidth]{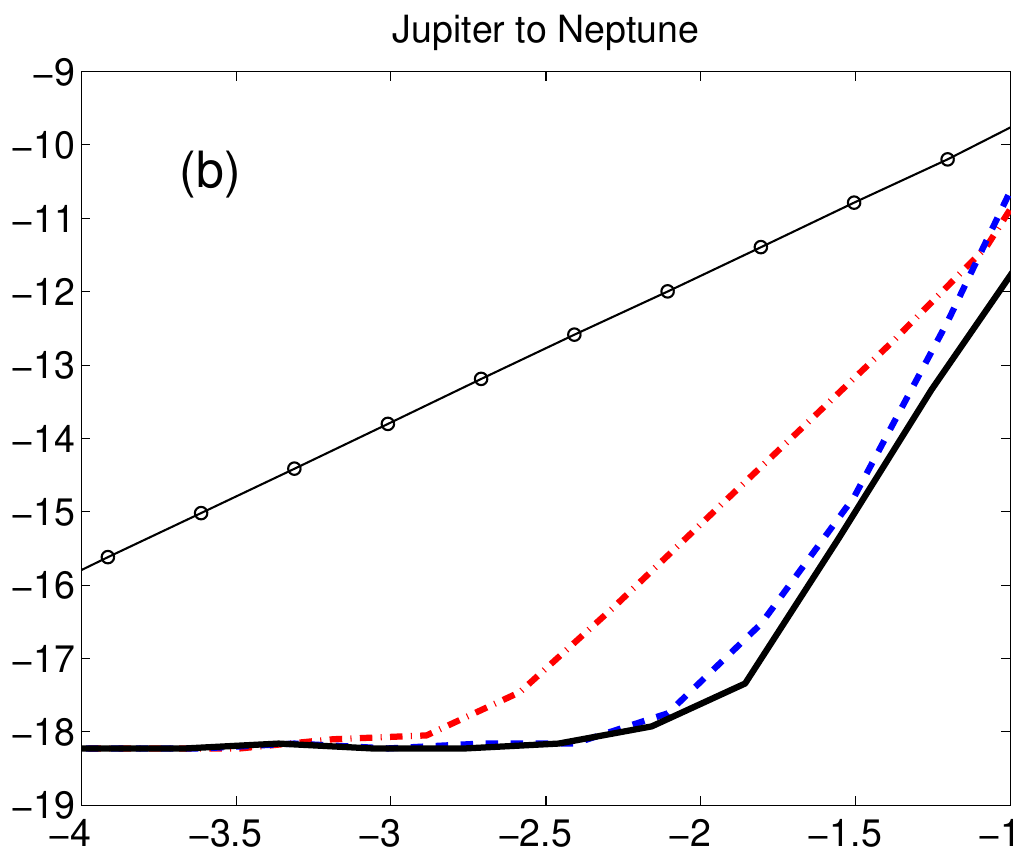}
\includegraphics[width = 0.49\textwidth]{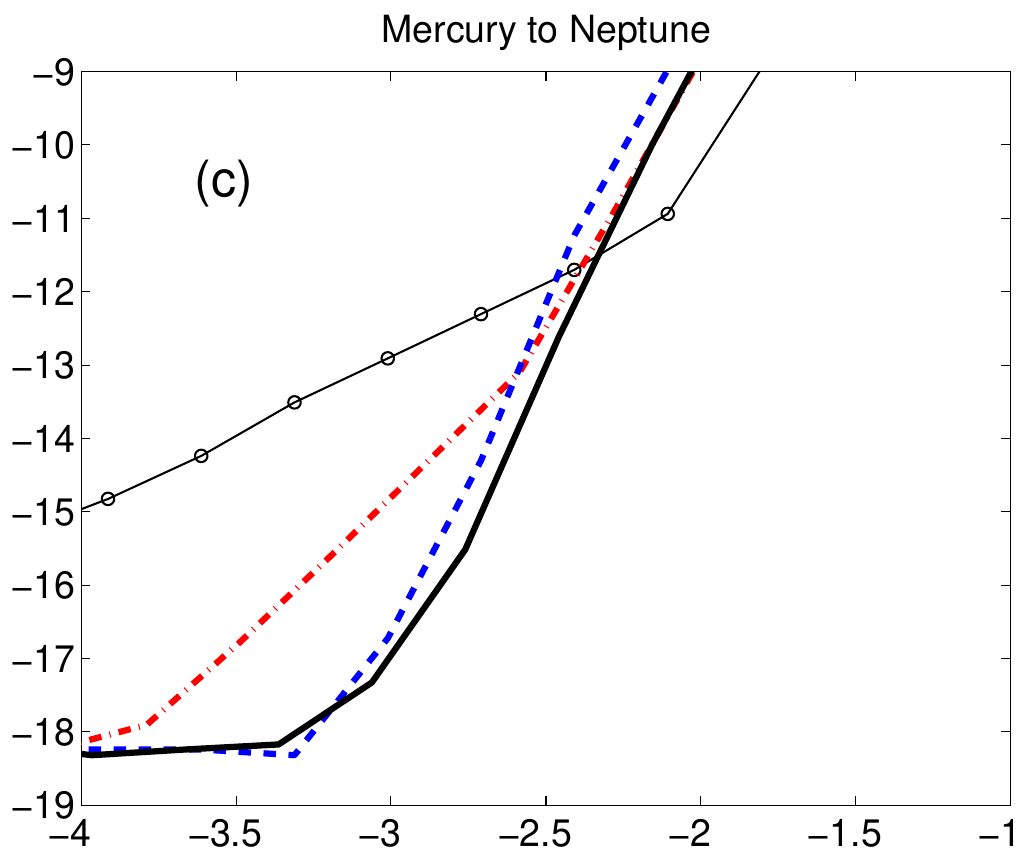}
\caption{ \small{Comparison between ABA schemes of order (8,4), (8,6,4)
and (10,6,4) of the form (\ref{eq:3approx}), with the approximate flow 
$\widetilde{\varphi}_{\tau}^{[b]}$ given by the leapfrog composition
(\ref{helio.4}) and
\texttt{ABA82} (also with the approximation (\ref{helio.4})). Panel (a): the 4 inner planets; panel (b): the 4
outer planets and panel (c): the whole Solar System. The $x$-axis represents
the (inverse of the) cost  $\tau/s$, and the $y$-axis the maximum
energy variation for one integration with constant step size
$\tau$. Both are in logarithmic scale.} }\label{Fig_ABAH_high}
\end{figure}

In Figure~\ref{Fig_ABAH_high} we can see the results achieved by the new 
ABA splitting schemes of Table~\ref{TAB_ABAH_high} on the N-body problem 
in Heliocentric coordinates for different planetary configurations. Method \texttt{ABA82} refers here to the
composition (\ref{eq:3approx}) with the coefficients of the ABA scheme of generalized order $(8,2)$ given in 
\cite{laskar01hos,mclachlan95cmi}, but with the flow $\varphi_{\tau}^{[b]}$
approximated by the leapfrog (\ref{helio.4}).
Diagram (a) corresponds to the four inner planets (Mercury, Venus, Earth and 
Mars), picture (b) to the four outer planets (Jupiter, Saturn, Uranus and Neptune) 
and finally diagram (c) is obtained when the eight planets in 
the Solar System (Mercury to Neptune) are taken into account. Initial conditions and mass parameters have been taken from the INPOP10a 
planetary ephemerides \cite{2011CeMDA.111..363F}
 (\verb+http://www.imcce.fr/inpop/+). 

We have integrated the same initial conditions for each scheme 
using different step sizes $\tau$. For each step size
($\tau_i = 1/2^i$ for $i = 1, 15$) we have computed the numerical trajectory over 
{\tt niter} = $10^5$ evaluations of the integration scheme (i.e. if $\tau = 0.5$, then the final time is 
$t_f = 50000$years). For each trajectory we plot the maximum variation in energy along 
the trajectory versus the inverse of the computational cost, $\tau/s$, both in logarithmic 
scale. 

All the simulations have been done in Fortran using extended double precision and 
compensated summation during the evaluation of the inner
stages of each scheme.

Notice that in the case of the four inner planets (Figure~\ref{Fig_ABAH_high}  (a))
the performance for the different ABA schemes in Table~\ref{TAB_ABAH_high} are better 
than method \texttt{ABA82}, but there is not much 
difference between the ABA schemes \texttt{ABAH844}, \texttt{ABAH864} and \texttt{ABAH1064}.
Nevertheless, 
if we look at the results for the four outer planets and the whole Solar System 
(Figure~\ref{Fig_ABAH_high} (b) and (c) respectively) we can see that 
\texttt{ABAH864} and \texttt{ABAH1064} show better results than \texttt{ABAH844}.  

Despite the fact that the size of the perturbation for the different planetary configurations 
presented can be very different \cite{farres12sif}, we believe that the difference between the 
performance of the schemes for the inner and the outer planets in the Solar System is 
mainly due to Mercury. Its fast orbital period and relatively high eccentricity are the main limiting 
factors when one tries  to improve the efficiency of the higher order schemes. Notice that 
for the inner planets the orbital period is much shorter than for the outer planets, and so,
according with \cite{viswanath02hmt}, this imposes a restriction on the step size to be used
and the number of evaluations per orbital period.

\section{Concluding remarks}

In reference \cite{laskar01hos}, symplectic splitting methods of generalized order $(2n,2)$ 
up to $n=5$ (first described in \cite{mclachlan95cmi}) were
systematically derived and tested on the Sun--Jupiter--Saturn system over 25000 years in Jacobi coordinates,
observing an improvement in the accuracy with respect to the leapfrog integrator by several
orders of magnitude at the same computational cost. Methods in this family have all the coefficients positive and good
stability properties. Scheme $(8,2)$ in particular has been used in several long
term simulations of the whole Solar System (e.g., \cite{laskar09eoc,laskar04alt}) and corresponds to method
\texttt{ABA82} in the examples reported here. All the tests carried out in \cite{laskar01hos} showed that the error
term $\eps^2 \tau^3$ was the main limiting factor in the performance of the integrators, so the natural
question was whether schemes of higher order (and thus already involving some negative coefficients) could
be useful for integrating planetary N-body problems.

As a matter of fact, methods of order $(8,4)$ obtained in \cite{mclachlan95cmi} do improve the performance of
\texttt{ABA82} for this problem in Jacobi coordinates, as the experiments reported
in \cite{farres12sif} show. In this work we have pursued this line of research and constructed new families
of higher order splitting methods specifically oriented to the numerical integration of near-integrable Hamiltonian
systems, and in particular for planetary N-body problems, both in Jacobi and Poincar\'e Heliocentric coordinates.
For this purpose, first we have derived explicitly the set of independent necessary and sufficient order conditions that
splitting methods must verify to achieve a certain order of accuracy and then we have solved these equations. A non-trivial
task that requires the use of homotopy continuation techniques and optimization criteria to select the most appropriate
solution.

Although the new methods involve some negative coefficients, and thus one could think that their numerical
stability might be compromised, they have been selected to minimize the error terms at higher orders and
the sum of the absolute values of their coefficients. As a result, the size of the negative coefficients of our new methods is
relatively small. In any case, the experiments reported here clearly indicate that the new methods of order
$(8,6,4)$ and $(10,6,4)$ achieve accuracy up to round off error with larger step sizes than 2nd-order schemes.

There are near-integrable systems of the form (\ref{eq:1}) where the exact flow $\varphi_{\tau}^{[b]}$
corresponding to the perturbation
is not available. In that case, we have constructed splitting methods of the form (\ref{eq:3approx}), where
$\widetilde{\varphi}_{\tau}^{[b]}$ is a 2nd-order symmetric approximation of the exact flow. This class of schemes
has shown to be particularly efficient for long time integrations of N-body planetary systems in Poincar\'e Heliocentric coordinates
when the leapfrog approximation (\ref{helio.4}) is considered.

Numerical simulations show that the efficiency of the new integrators presented here is essentially similar in both 
Jacobi and Heliocentric coordinates.
We believe this result is worth remarking, since canonical Heliocentric variables provide very often a more
convenient formulation of the problem. The improvement of the new integrators presented here (in particular, methods
of order $(8,6,4)$ and $(10,6,4)$) with respect to
previous schemes is most notably exhibited when they are applied for the numerical integration of 
the outer planets. When the whole Solar System is considered, although methods
of order $(8,6,4)$ and $(10,6,4)$ still provide the best results, it
is Mercury with its relatively high eccentricity and fast orbital
period  which constitutes the main limiting factor in all simulations.

When designing the new methods of generalized order $(8,6,4)$ and $(10,6,4)$, 
we have only considered compositions with the minimum number of stages to solve all
the order conditions. It might be the case, as in other contexts, that introducing more stages with additional free parameters 
could lead to more efficient schemes. We intend to explore this possibility and eventually collect the new
methods obtained, both in the ABA and BAB classes, in our website (\texttt{www.gicas.uji.es/software.html}).

\subsection*{Acknowledgements}

The work of SB, FC, JM and AM has been partially supported by Ministerio de
Ciencia e Innovaci\'on (Spain) under project MTM2010-18246-C03
(co-financed by FEDER Funds of the European Union), whereas AF and JL acknowledge financial support
by the FP7 GTSnext project.

\begin{appendix} 

\section*{Appendix}

To illustrate the numerical procedure we have followed to
obtain the methods with $s > 7$ stages in this work
we now describe in
detail the construction of method \texttt{ABAH1064}  of
generalized order (10,6,4), whose coefficients are collected in Table \ref{TAB_ABAH_high}. It has the form 
(\ref{eq:3approx}), 
where $\widetilde{\varphi}^{[b]}_{\tau}$ represents one step of some 2nd-order symmetric method used to approximate the
flow $\varphi^{[b]}_{\tau}$. This is the case, for instance, when Heliocentric coordinates are used for the integration
of the Solar System.

We consider nine-stage methods with the following sequence of coefficients:
\begin{equation}   \label{10.6.4b}
    a_1 \, b_1 \, a_2 \, b_2 \,  a_3 \, b_3 \, a_4 \, b_4 \, a_5 \, b_5 \, a_5 \, b_4 \, a_4 \, b_3 \,  a_3 \, b_2 \, a_2 \, b_1 \, a_1. 
\end{equation}
Symmetric splitting methods of generalized order $(10,6,4)$ for
Heliocentric coordinates must satisfy ten order conditions, that is, consistency 
\[
a_1 + a_2 +a_3 +a_4 + a_5 =\frac{1}{2}, \qquad
2(b_1 + b_2 +b_3 +b_4) + b_5 =1, 
\]
the special constraint (\ref{helio.5}) for Heliocentric coordinates, 
and the order conditions related to the Lyndon multi-indices (3), (5), (7), (9), (1,2), (1,4), (2,3) in
Table~\ref{tab.1} (recall that the $c_i$ are given by
(\ref{eq:ci})). We thus have ten polynomial equations and ten
unknowns, which we collect in a vector
$x=(a_1,\ldots,a_5,b_1,\ldots,b_5) \in \mathbb{R}^{10}$. The system of
algebraic equations one aims to solve can then be written in the
compact form $f(x)=0$.
Recall that any solution of such system must have at least
one negative $a_i$ and one negative $b_j$. We are interested in
finding solutions  with a small Euclidean norm
$\|x\|=\|(a_1,\ldots,a_5,b_1,\ldots,b_5)\|$ (to ensure that the negative $a_i$ and $b_j$ have small absolute
values).  

In order to do that, we first  split the system $f(x)=0$ into 
\begin{equation}    \label{eq:f1f2}
f_1(x)=0, \qquad\quad   f_2(x)=0,
\end{equation}
where $f_2(x)=0$ corresponds to the conditions for
  the Lyndon  multi-indices (1,2), (1,4), and (2,3), whereas $f_1(x)=0$
  collects the remaining seven equations. We
  then proceed as follows:

\begin{itemize}
\item We determine the point $x^0=(a^0_1,\ldots,a^0_5,b^0_1,\ldots,b^0_5)  \in
\mathbb{R}^{10}$ as the (unique) solution of the following constrained minimization problem:
\[
\min_{a_3=a_4=0,  f_1(a_1,\ldots,a_5,b_1,\ldots,b_5)=0} \; \sum_{i=1}^{5}(b_i^2+a_i^2).
\] 
It is not difficult to check that the sequence of coefficients 
\[
   a_1^0 \, b_1^0 \, a_2^0 \, (b_2^0+b_3^0+b_4^0) \,  a_5^0 \, b_5^0
   \, a_5^0 \, (b_2^0+b_3^0+b_4^0)\, a_2^0\, b_1^0\, a_1^0
\]
corresponds precisely to  the 5-stage symmetric ABA
  method of generalized order (10,2) with positive coefficients
  considered in \cite{laskar01hos,mclachlan95cmi,mclachlan96mos}. Thus, it is our starting point in the search
  of an efficient (10,6,4) method.

\item We then choose an arbitrary orthogonal matrix $M \in
  \mathbb{R}^{3\times 10}$, and for a randomly chosen complex number
  $\gamma \in \mathbb{C}$, consider the following one-parameter family
  of systems of polynomial equations:
\begin{equation}
  \label{eq:f1f2b}
  f_1(x)=0, \qquad\quad t\, f_2(x) + (1-t)\, \gamma\, M \cdot (x-x^0) = 0.
\end{equation}
For a generic $M$ and $\gamma$, there exists a unique continuous curve
$x=\rho(t) \in \mathbb{C}^{10}$ of solutions of this family
of polynomial systems  such that $\rho(0)=x^0$. Here $t \in [0,1)$ denotes the continuation parameter. If 
$x=\rho(1) \in \mathbb{R}^{10}$, then $x$ is a real
solution of the original system (\ref{eq:f1f2}). We have applied a
numerical continuation algorithm to compute such a solution for
several values of $\gamma \in \mathbb{C}$, and found two real
solutions. The solution $x$ with smaller norm $\|x\|$
gives the method \texttt{ABAH1064}  of
generalized order (10,6,4) displayed in Table~\ref{TAB_ABAH_high}. Notice the small absolute value
of both $a_3$ and $a_4$ in the resulting scheme (recall that this solution has been obtained starting with
$x^0$ such that $a_3^0=a_4^0=0$).
\end{itemize}

We have followed a similar procedure to compute solutions starting with
$x^0$ such that $a_i^0=a_j^0=0$ for indices $(i,j) \neq (3,4)$. 
Such procedure leads to non-real solutions $x$ for
some of the choices of $(i,j)$, and for
other choices gives real solutions
with larger norm $\|x\|$ and larger error terms than  method \texttt{ABAH1064}.

\end{appendix}

\bibliographystyle{plain}

\end{document}